\documentclass[a4paper,12pt]{article}

\usepackage[utf8]{inputenc}
\usepackage{amsthm, amssymb, mathrsfs, amscd, amsmath, stmaryrd}
\usepackage{mathtools}
\usepackage[pdfa]{hyperref}
\usepackage{graphicx}
\usepackage{indentfirst} %per avere il rientro anche dopo i titoli
\usepackage{multirow}
\usepackage{longtable}
\usepackage{amsfonts}
\usepackage{mathrsfs}
\usepackage{enumerate}
\usepackage{cite} %BibTeX
\usepackage{soul}
\usepackage{authblk}

\DeclarePairedDelimiter\absval{\lvert}{\rvert}
\newcommand{\abs}[1]{\absval*{#1}}

\newcommand{\sep}{\,\,|\,\,}

\theoremstyle{plain}
\newtheorem{Th}{Theorem}
\newtheorem{Prop}[Th]{Proposition}
\newtheorem{Lem}[Th]{Lemma}
\newtheorem{Cor}[Th]{Corollary}

\theoremstyle{definition}
\newtheorem{Def}[Th]{Definition}

\title{Vanishing ideals of binary Hamming spheres}

\author[1]{Giuseppe D'Alconzo}
\author[2]{Alessio Meneghetti}
\affil[1]{Department of Mathematical Sciences, Politecnico di Torino}
\affil[2]{Department of Mathematics, Università di Trento}
\date{}                     %% if you don't need date to appear
\begin{document}
\maketitle

\begin{abstract}
	We consider the simplified Algebraic Normal Form (sANF) of Boolean functions vanishing on Hamming spheres centred at zero and the associated sANF vector. We show that this vector is periodic, leading to an efficient computation of the sANF and to specific formulas for particular cases. Moreover, we explicitly provide a connection to the binary M{\"o}bius transform of the elementary symmetric functions. We conclude by presenting a method based on polynomial evaluation to bound the minimum distance of binary nonlinear codes. The same method can be used to compute the minimum distance and the weight distribution of binary linear codes.
\end{abstract}

\section{Introduction}

% 	Many computationally hard problems can be described by Boolean polynomial systems, and the standard approach is the computation of the Gr{\"o}bner Basis of the corresponding ideal. Since it is a quite common scenario, we will restrict ourselves to ideals of $\mathbb{F}_2[x_1, \dots,x_n]$ containing the entire set of field equations $\{x^2_i + x_i \}_i$. To ease the notation, our work environment will therefore be the quotient ring $R = \mathbb{F}_2[x_1, \dots,x_n]/(x^2_1 +x_1 , \dots, x^2_n +x_n )$. Moreover, most of our results do not depend on the number $n$ of variables, and when not otherwise specified we consider $R$ to be defined in infinitely many variables. We denote with $X$ the set of our variables.

In this work we characterise the vanishing ideal $I_t$ of the set of binary vectors contained in the Hamming sphere of radius $t - 1$. This characterisation corresponds to the explicit construction of the square-free polynomial $\phi_t$ whose roots are exactly the set of points of Hamming weight at most $t-1$. It is worth mentioning that this polynomial corresponds to the Algebraic Normal Form (ANF) of the Boolean function which vanishes inside the Hamming sphere. 
A direct application of our work is the possibility to add $\phi_t$ to the generating system of an ideal, and therefore to force the corresponding variety to live in the Hamming sphere of radius $t - 1$. For example, consider a linear code $C$, its parity-check matrix $H$, an integer $t$ and a syndrome $s$. The Maximum Likelihood Problem (MLD) \cite{berlekamp1978inherent} is the NP-complete problem of determining whether a given vector is within Hamming distance $t$ from a codeword of a given code, or, more precisely, whether there exists a vector $v$ of Hamming weight less than $t$ such that $Hv^{\top}=s^{\top}$. Thus, MLD can be described as the problem of deciding whether the variety associated to the system
$$
\left\{
\begin{array}{l}
	Hv^{\top}=s^{\top}\\
	\phi_t(v)=0
\end{array}
\right.
$$
is non-empty. We point out that this description of MLD does not readily provide a useful formulation of the problem in terms of polynomial ideals, since $\phi_t$ is a dense polynomial. An interested reader may refer to \cite{meneghetti2021explicit} for a description on MLD in terms of multivariate quadratic polynomial systems.
A less straightforward application will be presented in Section \ref{Sect:Applic}, where we show a novel method to check whether the minimum distance of a linear code is in a given range. We believe that the ideas presented in Section \ref{Sect:Applic} could eventually be a starting point for new algebraic algorithms for the computation of the minimum distance of linear codes. This would however require the design of dedicated procedures to minimise the computational complexity of such algorithms, and this is beyond the aim of this paper. The reader can find similar methods in \cite{guerrini2010computing}, where the authors proposed a technique to compute the distance distribution of systematic non-linear codes by relying on polynomial ideals. Other interesting results, obtained to deal with particular classes of codes, can be find in \cite{garcia1999algorithm} and \cite{hu2004computation}. The main difference between previously known algorithms and the ideas presented in this work is that we do not need to rely on brute-force-like methods, nor we require the computation of a Gr{\"o}bner basis. An interested reader can find in \cite{carlet2010boolean} a thorough discussion about Boolean functions and in \cite{canteaut2005symmetric} an analysis of symmetric Boolean functions and their properties. For a comprehensive work on the utilisation of computational methods to address problems in algebra and geometry see \cite{cox1994ideals}, while for everything regarding Coding Theory we refer to \cite{macwilliams1977theory}.

The paper is organised as follows.
In Section \ref{Sect:BinarySymmFunc} we introduce the notation and provide some properties of binary symmetric functions, and in Section \ref{Sect:Moebius} we discuss the binary M{\"o}bius transform. These preliminary results will then be used in the remaining Sections \ref{Sect:VanishingId} and \ref{Sect:Applic}, where we respectively discuss about the generating polynomials $\phi_t$ of the vanishing ideals of Hamming spheres, and their application to the computation of the distance of linear codes. The \hyperref[Sect:Appendix]{Appendix} analyses the structure of $\phi_t$ for some particular forms of $t$.

% 	Many computationally hard problems can be described by Boolean polynomial systems, and the standard approach is the computation of the Gr{\"o}bner Basis of the corresponding ideal. 
\subsection{Some remarks on our notation}
Let $(\mathbb{F}_2,+,\cdot)$ be the finite field with two elements. In this work we only work with ideals of $\mathbb{F}_2[x_1, \dots,x_n]$ containing the entire set of field equations $\{x^2_i + x_i \}_i$. To ease the notation, our work environment will therefore be the quotient ring $R = \mathbb{F}_2[x_1, \dots,x_n]/(x^2_1 +x_1 , \dots, x^2_n +x_n )$, and thus we do not list the fields equations among the generator polynomials of ideals. Similarly, when we work with the varieties associated to ideals we omit to explicitly state that we restrict them to $\mathbb{F}_2$. Moreover, most of our results do not depend on the number $n$ of variables, and when not otherwise specified we consider $R$ to be defined in infinitely many variables.
\\
In some results and proofs we consider operations over the integers, e.g. when dealing with indices or binomial coefficients, using with a slight abuse of notation the symbols $+$ and $\cdot$. A sum is therefore either the sum of $\mathbb{F}_2$ or $\mathbb{Z}$ according to the context.
\\
Following the notation adopted in Coding Theory we consider any vector $v$ to be a row vector, and when we are required to work with column vectors we denote them as $v^{\top}$. The Hamming weight of a vector $c$ is denoted with $\mathrm{w}(c)$ and is the number of its non-zero coordinates. The Hamming distance between two vectors $b$ and $c$ is the Hamming weight of $b+c$ and is denoted with $\mathrm{d}(b,c)$.

\section{Binary symmetric functions}\label{Sect:BinarySymmFunc}

The vanishing ideal of a Hamming sphere of radius $t-1$ centred at zero is generated in $R = \mathbb{F}_2[x_1, \dots,x_n]/(x^2_1 +x_1 , \dots, x^2_n +x_n )$ by a single binary polynomial which we denote with $\phi_t$. To be precise, $\phi_t$ is the Algebraic Normal Form of the map $\left(\mathbb{F}_2\right)^n \to \mathbb{F}_2$ whose zeros are all and only the binary vectors $v$ whose Hamming weight $\mathrm{w}(v)$ is less than $t$, i.e. $v$ has less than $t$ ones. By definition, $\phi_t$ depends uniquely on the weight of the input, hence $\phi_t$ is a symmetric polynomial, which can therefore be written in terms of the elementary symmetric functions (ESFs) $\sigma_i = \sum_{m\in M_i}m$, where $M_i$ is the set of all monomials of degree $i$. We recall that in our set-up each $m\in M_i$ is a square-free monomial; this implies that any monomial of degree $i$ is the multiplication of $i$ distinct variables. In case of $n$ variables, $\abs{M_i} = \binom{n}{i}$.

In this section we look more closely to the behaviour of ESFs when working in the quotient $R$. The properties described in this section will then be used in Section \ref{Sect:VanishingId} to derive further results on $\phi_t$.

Let $a=(a_0,\dots, a_k)$ and $b=(b_0,\dots, b_k)$ be two bit-strings of equal length (otherwise we can pad the shortest with 0's). We define an order $\preceq$ on bit-strings, where $a\preceq b$ if and only if for every $i$ in $\{0,\dots, k\}$ we have that $a_i\le b_i$ (i.e. on single bits we use the order induced by integers: $0\le 0,\; 0\le 1,\; 1\le 1$ and $1>0$).

\begin{Def}
	We denote with $\mathrm{b}(n)$ the binary representation of the non-negative integer $n$.
\end{Def}

With an abuse of notation we write $n_1\preceq n_2$ with $n_1$ and $n_2$ non-negative integer when $\mathrm{b}(n_1)\preceq\mathrm{b}(n_2)$.

The following is a well known result on Boolean functions. \cite{carlet2010boolean}.
\begin{Th}
	Let $f$ be a Boolean function in $n$ variables. Then there exists a Boolean function $g$ in $n$ variables such that
	$$f(x_1,\dots,x_n)=\sum_{(a_1,\dots,a_n)\in \left(\mathbb{F}_2\right)^n} g(a)x_1^{a_1}\cdots x_n^{a_n}.$$
	This writing is called \emph{Algebraic Normal Form} (ANF) of $f$, the function $g$ is called \emph{binary M{\"o}bius Transform} of $f$ and is denoted with $g=\mu(f)$.
\end{Th}

It is known \cite{carlet2010boolean} that the binary M{\"o}bius transform has order 2: for every $f$ we have $\mu(\mu(f))=f$.

\begin{Def}
	The elementary symmetric function of degree $i$ in $n$ variables is defined as
	$$\sigma_{i} = \sum_{m\in M_i} m.$$
\end{Def}

Since every symmetric function is a polynomial in the variables $\sigma_{i}$, we have that for $f$ in $R$ we can write $f$ as a linear combination of elementary symmetric functions $\{\sigma_{i}\}_i$.
\begin{Def}\label{sANF}
	Every symmetric Boolean function $f$ in $R$ can be written in the following form called \emph{simplified Algebraic Normal Form} (sANF):
	$$f = \sum_{i=0}^n \lambda_{f}(i) \sigma_{i}$$
	where the vector $\lambda_{f}=(\lambda_{f}(0), \dots, \lambda_{f}(n))$ in $\left(\mathbb{F}_2\right)^{n+1}$ is called sANF vector of $f$.
\end{Def}

The sANF of a symmetric Boolean function $f$ is a useful tool since it allows to reduce the length of the representation of a function from $2^n$ to $n+1$ bits.

\begin{Th}[Lucas' Theorem]\label{Th:Lucas}
	Let $i$ and $w$ be two non-negative integers. Then
	$$\binom{w}{i} = 1 \mod 2 \iff \mathrm{b}(i)\preceq \mathrm{b}(w). $$
\end{Th}

\begin{Th}\label{ThProd}
	Let $\sigma_{i}$ and $\sigma_{j}$ be defined in the quotient $R$. Then $\sigma_{i}\cdot\sigma_{j}=\sigma_h$ with $\mathrm{b}(h)=\mathrm{b}(i)\vee\mathrm{b}(j)$.
\end{Th}
\begin{proof}
	Let $i < j$. From Definition \ref{sANF} it follows that the product of $\sigma_{i}$ and $\sigma_{j}$ is a linear combination of ESFs. Given $m_1$ and $m_2$ in $R$, the degree of $m_1m_2$ is at least $\max(\deg(m_1), \deg(m_2))$ and at most $\deg(m_1) + \deg(m_2 )$, namely $\sigma_{i}\sigma_{j} = \sum_{h=j}^{i+j} c_h \sigma_{h}$. Since $\sigma_{i}$ and $\sigma_{j}$ are symmetric polynomials, if a monomial $m_h$ of degree $h$ appears in their product, then $c_h = 1$. Observe that $m_h$ is a monomial in $\sigma_{i}\sigma_{j}$ if and only if the cardinality of the set $\{(m_i , m_j ) : m_im_j = m_h , m_i \in M_i , m_j \in M_j \}$ is odd. The number of these pairs $(m_i, m_j)$ is equal to
	\begin{equation}\label{Eq1}
		N = \binom{h}{i}\binom{i}{i+j-h},
	\end{equation} 
	where the first term is equal to the ways of choosing $i$ variables (the monomial $m_i$) from a set of $h$ variables (the monomial $m_h$), and the second term corresponds to choosing $j - (h - i)$ variables among the $i$ which appear in $m_i$ (obtaining the monomial $m_j$).
	
	Since $\binom{n}{k} = \binom{n}{n-k}$, the second binomial coefficient in Equation \eqref{Eq1} can also be written as $\binom{i}{h-j}$. The product $N$ in Equation \eqref{Eq1} is odd if and only if both binomial coefficients are odd, hence, by Theorem \ref{Th:Lucas}, if $\mathrm{b}(i) \preceq \mathrm{b}(h)$ and $\mathrm{b}(h-j) \preceq \mathrm{b}(i)$. We write $\mathrm{b}(h)$ as $\mathrm{b}(i) \vee \mathrm{b}(j) + b_1 + b_2$, with $b_1, b_2$ binary vectors such that $b_1 \preceq \mathrm{b}(i) \vee \mathrm{b}(j)$ and $b_2 \vee (\mathrm{b}(i) \vee \mathrm{b}(j)) = 0$. We have three possible cases:
	\begin{itemize}
		\item[1.] $b_1=b_2=0$;
		\item[2.] $b_1=0$ and $b_2\ne0$;
		\item[3.] $b_1\ne 0$.
	\end{itemize}
	In the first case, $\mathrm{b}(h) = \mathrm{b}(i) \vee \mathrm{b}(j)$, then both conditions are satisfied and this implies that $c_h = 1$.\\
	In the other cases, at least one of the binomial coefficients in Equation \eqref{Eq1} is even. This implies that if $\mathrm{b}(h) \ne \mathrm{b}(i)) \vee \mathrm{b}(j)$ then $c_h = 0$.
\end{proof}

While Definition \ref{sANF} assures us that any binary symmetric polynomial in $R$ can be written as a linear combination of $\{1, \sigma_{1},\dots, \sigma_{n}\}$, an interesting consequence of Theorem \ref{ThProd} is that to represent the same polynomial we do not really need all the ESFs. All such polynomials can indeed be defined in terms of sums and products of elements in $\{1, \sigma_{1} , \sigma_{2} , \sigma_{2^2} , \dots, \sigma_{2^s} ,\dots\}$.

\begin{Cor}
	The set of all square-free binary symmetric polynomials in $n$ variables is equivalent to $\mathbb{F}_2[y_0,\dots,y_s]$, with $s=\lfloor\log_2(n)\rfloor$.
\end{Cor}
\begin{proof}
	We use Definition \ref{sANF} and Theorem \ref{ThProd}. Equivalently, we can start from a polynomial $f$ in $\mathbb{F}_2[y_0,\dots,y_s]$ and compute $f(\sigma_{1},\sigma_{2},\dots,\sigma_{2^s})$.
\end{proof}

\section{The binary M{\"o}bius transform and the simplified value vector}\label{Sect:Moebius}

In this section we present the simplified value vector of a symmetric Boolean function and its relation with the sANF. Moreover we present a closed formula for the binary M{\"o}bius transform, useful to compute the sANF.

\begin{Th}
	Let $f:\left(\mathbb{F}_2\right)^n\to \mathbb{F}_2$ be a Boolean function. Then its binary M{\"o}bius transform is
	\begin{equation}\label{moeb}
		\mu(f)(x_1,\dots,x_n) = (1+x_1)\cdots (1+x_n)\cdot f\left( \frac{x_1}{1+x_1},\dots,\frac{x_n}{1+x_n} \right),
	\end{equation}
	where the fractions are symbolic, since their denominators vanish together with the corresponding term in the product on the left.
\end{Th}

We provide an example before proving the formula. Let $f(x_1 , x_2 ) = x_1 +x_1 x_2$. Then
\begin{align*}
	\mu(f)(x_1,x_2) &= (1+x_1)(1+x_2)\left[ \frac{x_1}{1+x_1} + \left(\frac{x_1}{1+x_1}\right)\left(\frac{x_2}{1+x_2}\right) \right] = \\
	% 		& = (1+x_1)(1+x_2)\frac{x_1}{1+x_1} + (1+x_1)(1+x_2)\left(\frac{x_1}{1+x_1}\right)\left(\frac{x_2}{1+x_2}\right)= \\
	& = (1+x_2)x_1 + x_1x_2,
\end{align*}
namely $\mu(f)=x_1$.

\begin{proof}
	The binary M{\"o}bius transform of $f$ is the Boolean function whose evaluation vector corresponds to the coefficients of the ANF of $f$. To be more precise, a point $b = (b_1 ,\dots, b_n ) \in (F_2 )^n$ can be identified with the monomial $m_b = X^b = x_1^{b_1}  \cdots x^{b_n}_n$, where $x^1_i$ = $x_i$, while with $x^0_i$ we mean that $x_i$ does not appear in $m_b$. Hence, by definition, $\mu(f)(b) = 1$ if and only if $m_b = X^b$ is a monomial of $f$.
	
	We consider then a generic monomial $m_b = X^b$ in $f$, and we observe that by formula \eqref{moeb} we obtain
	$\mu(m_b) = (1 + x_1 )^{1+b_1} \cdots (1 + x_n )^{1+b_n} \cdot x^{b_1}_1 \cdots x^{b_n}_n$, namely we obtain $\mu(m_b) = (1 + X)^{1+b} \cdot X^b$. It can easily be checked that the polynomial $\mu(m_b)$ assume value $1$ only when evaluated at $b$. Then, to each monomial $m_b$ in $f$ it corresponds a polynomial $\mu(m_b)$ in $\mu(f)$ for which $\mu(m_b)(b)=1$, while it is zero everywhere else. So the evaluation vector of $\mu(f)$ is exactly the vector of coefficients of $f$.
\end{proof}

If $f$ is a symmetric Boolean function, its values depend only on the Hamming weight of the input \cite{canteaut2005symmetric}:
\begin{Prop}
	Let $f$ be a symmetric Boolean function. Then there exists
	$$v_f : \{0,\dots, n\} \to \mathbb{F}_2$$
	such that $f(x)=v_f(\mathrm{w}(x))$.
\end{Prop}
From this we can define the truth table of a symmetric function.
\begin{Def}
	The \emph{simplified value vector} of a symmetric Boolean function $f$ is denoted with
	$$v_f = (v_f(0),\dots, v_f(n))\in \left(\mathbb{F}_2\right)^{n+1}.$$ 
\end{Def}

\begin{Prop}\label{Propsum}
	Let $f$ and $g$ be two symmetric Boolean functions in $n$ variables. Then
	\begin{itemize}
		\item[(i)] $v_{f+ g} = v_f + v_g$, and 
		\item[(ii)] $\lambda_{f+ g} = \lambda_f + \lambda_g$.
		%\item[2)] $\lambda_{fg} = \lambda_f \cdot \lambda_g$, where $\cdot$ denotes the pairwise product.
	\end{itemize}
\end{Prop}
\begin{proof}
	Let $f$ and $g$ be two symmetric Boolean functions in $n$ variables.\\
	The equation $(i)$ follows from $(f+ g)(c)=f(c) + g(c) = v_f(i) + v_g(i)$, where $c$ is a vector with $\mathrm{w}(c)=i$.\\
	Result $(ii)$ can be proven writing
	$$f+g = \sum_{i=0}^n \lambda_{f+g}\sigma_i $$
	and using the definition of sANF vectors $\lambda_f$ and $\lambda_g$ of $f$ and $g$, respectively.
\end{proof}

A link between the sANF vector and the simplified value vector of a symmetric Boolean function $f$ is given in the following result.

\begin{Prop}[{\cite[Proposition~2]{canteaut2005symmetric}}]\label{PropsANF}
	Let $f$ be a symmetric Boolean function with sANF vector given by $(\lambda_{f}(0), \dots, \lambda_{f}(n))$ and simplified value vector $(v_f(0),\dots, v_f(n))$. Then for every $i$ in $\{0,\dots, n\}$:
	\begin{itemize}
		\item[1.] $v_f(i) = \sum_{k\preceq i} \lambda_{f}(k)$;
		\item[2.] $\lambda_f(i) = \sum_{k\preceq i} v_{f}(k)$.
	\end{itemize}
\end{Prop}

The previous proposition give us a way to compute one of the two characteristic vectors of a Boolean function when the other is known.

\begin{Lem}\label{LemMob}
	Let $f$ be a symmetric Boolean functions and $\mu(f)$ its binary M{\"o}bius transform. Then $v_f=\lambda_{\mu(f)}$ and $\lambda_f=v_{\mu(f)}$.
\end{Lem}
\begin{proof}
	In the case of arbitrary Boolean functions, the binary M{\"o}bius transform $g=\mu(f)$ of $f$ satisfies the following relation in the ANF of $f$:
	$$ f(x_1,\dots, x_n)=\sum_{(a_1,\dots,a_n)\in \left(\mathbb{F}_2\right)^n} g(a)x_1^{a_1}\cdots x_n^{a_n}.$$
	If $f$ is symmetric, also its binary M{\"o}bius transform is symmetric and the sANF of $f$ gives us that the simplified value vector of $\mu(f)$ is equal to $\lambda_{f}$, leading to $v_{\mu(f)}=\lambda_{f}$.\\
	The same argument can be applied to $g=\mu(f)$. In fact, using $\mu(\mu(f))=f$, we have $v_{f}=\lambda_{\mu(f)}$.
\end{proof}
\iffalse	
Finally we study the behaviour of the simplified value vector and the binary M{\"o}bius transform when more Boolean functions are composed.
\begin{Prop}
	Let $f:\left(\mathbb{F}_2\right)^n\to \mathbb{F}_2$ and let $g_1,\dots,g_n$ be Boolean functions such that $g_i:\left(\mathbb{F}_2\right)^k\to \mathbb{F}_2$ for every $i$ in $\{1, \dots , n\}$. Consider the composition
	$$ 
	\begin{array}{rccl}
		F:&\left(\mathbb{F}_2\right)^k&\to&\mathbb{F}_2\\
		&(x_1,\dots,x_k)&\mapsto & f(g_1(x_1,\dots, x_k),\dots, g_n(x_1,\dots,x_k)),
	\end{array}
	$$
	then
	$$F(u) = \sum_{w \preceq \left(g_1(u),\dots,g_n(u)\right)} \lambda_f(w)$$
	and $\lambda_F(u) = \sum_{a\preceq u} F(u)$.
\end{Prop}
\begin{proof}
	For each $u$ in $\left(\mathbb{F}_2\right)^k$ we have
	\begin{align*}
		F(u) &= f(g_1(u),\dots,g_n(u)) = \\
		& = \sum_{w\in \left(\mathbb{F}_2\right)^n} \lambda_f(w)g_1(u)^{w_1}\cdots g_n(u)^{w_n} 
	\end{align*}
	where in the last equality we used the ANF of $f$. Now we can observe that all the terms having at least one $i\in \{1,\dots,n\}$ such that $g_i(v)=0$ and $w_i=1$ vanish. Therefore we can write
	$$F(u) = \sum_{w \preceq \left(g_1(u),\dots,g_n(u)\right)} \lambda_f(w).$$
	It is known \cite{carlet2010boolean} that, for every Boolean function $h$, we have $\lambda_h(u)=\sum_{a\preceq u} h(a)$, hence we obtain $\lambda_F(u) = \sum_{a\preceq u} F(a)$.
\end{proof}
\fi

\section{The Vanishing ideal of a Hamming sphere}\label{Sect:VanishingId}

In this section we provide a description of $\phi_t$ and we prove some periodicity properties. Definition \ref{sANF} allows us to write
$$ \phi_t = \sum_{i=0}^n \lambda_{\phi_t}(i)\sigma_{i},$$
and we want to compute explicitly the coefficients of the sANF vector $\lambda_{\phi_t}$.
Proposition \ref{PropsANF} allows to write a formula for $\phi_t$ directly from its definition: since its simplified value vector has $t$ 0's at the beginning and then $n-t$ 1's then we can directly compute the sANF of $\phi_t$. Moreover we have the following results:

\begin{Lem}\label{Lems}
	Let $\lambda_{\phi_t}$ be the sANF vector of $\phi_t$. Then
	\begin{itemize}
		\item[1.] $\lambda_{\phi_t}(i)=0$ for each $i<t$;
		\item[2.] $\lambda_{\phi_t}(t)=1$;
		\item[3.] $\lambda_{\phi_t}(2^{\lceil\log_2t\rceil}) = 1$.
	\end{itemize}
\end{Lem}
\begin{proof}
	Since the simplified value vector of $\phi_t$ is $(0,\dots,0,1,\dots)$, with $t$ initial 0's, $1.$, $2.$ and $3.$ are direct consequence of Proposition \ref{PropsANF}.
\end{proof}

\begin{Th}
	Let $n\ge t$, and let $h_t=\sum_{i=t}^{n}\sigma_{i}$. Then $\phi_t = \mu(h_t)$.
\end{Th}
\begin{proof}
	We use Lemma \ref{LemMob}, observing that $v_{\phi_t}=(0,\dots,0,1,\dots,1)$ with $t$ initial $0$'s. Hence $\lambda_{\mu(\phi_t)} = v_{\phi_t}$ and $\mu(\phi_t)=h_t$. Applying the binary M{\"o}bius transform on both sides of the equation, we have $\phi_t=\mu(\mu(\phi_t))=\mu(h_t)$
\end{proof}

Let us study the periodicity of the sANF vector of $\phi_t$.

\begin{Def}
	Let $a$ be a vector in $\left(\mathbb{F}_2\right)^n$. We say that $a$ is \emph{periodic} with \emph{period} $T$ if for every $i$ in $\{1,\dots, T\}$ and every $k\ge 0$  such that $kT+i$ is less than $n$, we have that $a_i=a_{kT+i}$.
\end{Def}

In our case we need a slight generalised definition of a periodic vector.

\begin{Def}
	Let $a$ be a vector in $\left(\mathbb{F}_2\right)^n$. We say that $a$ is \emph{eventually periodic} with \emph{period} $T$ if there exists $\nu$ such that the vector $\tilde{a}$ in $\left(\mathbb{F}_2\right)^{n-\nu} $ defined by
	$$\tilde{a}_i=a_{i+\nu}\qquad i\in\{1,\dots,n-\nu\}$$ is periodic with period $T$. We say that the periodicity of $a$ starts from $\nu$.
\end{Def}

We observe that the two previous definitions can be extended to the case of infinitely long vectors, that, when considering the simplified value vector and the sANF vector of a Boolean function $f$, models the fact that $f$ has infinitely many variables.

\begin{Th}[{\cite[Theorem~1]{canteaut2005symmetric}}]\label{ThPer}
	Let $f$ be a symmetric Boolean function in $n$ variables with simplified value vector $v_f$. Then, $v_f$ is periodic with period $2^s<n$ if and only if $\mathrm{deg}(f)\le 2^s -1$. 
\end{Th}

Due to this result we can prove the following.

\begin{Th}\label{ThPerPhi}
	Let $s$ and $t$ such that $t$ is less than $n$ and $2^{s-1}\le t\le 2^s-1$. Then the sANF vector $\lambda_{\phi_t}$ of $\phi_t$ is eventually periodic with period $2^s$, starting from index $1$ and having $\lambda_{\phi_t}(0)=0$.
\end{Th}
\begin{proof}
	Consider the function $\phi_t+ 1$. Since it vanishes on vectors of weight at least $t$, its simplified value vector is given by
	$$v_{\phi_t+ 1} = (1,\dots, 1,0,\dots, 0)$$
	where there are $t$ initial $1$'s. Now let $\mu(\phi_t+ 1)$ be its binary M{\"o}bius transform, from Lemma \ref{LemMob} we have
	$$\lambda_{\mu(\phi_t+ 1)} = v_{\phi_t+ 1} = (1,\dots, 1,0,\dots, 0)$$
	and then we can write $\mu(\phi_t+ 1)$ as a polynomial
	$$\mu(\phi_t+ 1) = \sum_{i=0}^{t}\sigma_{i}.$$
	This means that the degree of $\mu(\phi_t+ 1))$ is $t\le 2^s-1$ and then Theorem \ref{ThPer} implies that its simplified value vector $v_{\mu(\phi_t+ 1)}$ is periodic of period $2^s$. This means, due to Lemma \ref{LemMob}, that $\lambda_{\phi_t+ 1}=v_{\mu(\phi_t+ 1)}$ is periodic of period $2^s$.\\
	By Proposition \ref{Propsum}, $\lambda_{\phi_t+ 1}=\lambda_{\phi_t}+ \lambda_1$ and we observe that $\lambda_1$ is the vector $(1,0,\dots, 0)$.
	This implies that $\lambda_{\phi_t}$ is eventually periodic of period $2^s$ starting from $1$. From Lemma \ref{Lems} we have $\lambda_{\phi_t}(0)=0$.
\end{proof}

The previous theorem shows that we need to compute only $2^{\lceil\log_2t\rceil}$ terms of the ANF of $\phi_t$, since they are periodic.\\
Due to the periodicity of the ANF of $\phi_t$ and the fact that it is independent from the number of variables $n$, we can always assume that a set of $n$ variables is obtained by restriction from a larger set. For this reason we can drop the notation of $n$.

Moreover, when $t$ is a power of $2$ we have the following result:
\begin{Cor}\label{Cor2s}
	Let $\lambda_{\phi_{2^s}}$ be the sANF vector of $\phi_{2^s}$, then $\lambda_{\phi_{2^s}}$ is eventually periodic of period $2^s$ starting from 1. Moreover
	$$\phi_{2^s} = \sigma_{2^s}+\sigma_{2\cdot2^s}+\sigma_{3\cdot2^s} + \dots = \sum_{j}\sigma_{j2^s}.$$
\end{Cor}
\begin{proof}
	From Theorem \ref{ThPerPhi} we have that $\lambda_{\phi_{2^s}}$ is eventually periodic of period $2^{s+1}$ starting from 1. From Lemma \ref{Lems} we have that for every $1<i<2^s$, $\lambda_{\phi_{2^s}}(i)=0$ holds. Now let $v_{\phi_{2^s}}$ be the simplified value vector of $\phi_{2^s}$, then
	$$\lambda_{\phi_{2^s}}(i+2^s) = \sum_{k\preceq i+2^s}v_{\phi_{2^s}}(k),$$
	since $v_{\phi_{2^s}}(k)=0$ for every $k<2^s$ we reduce the sum to $2^s\le k< i+2^{s}$ such that $k\preceq i+2^s$, and these are all equal to 1. Hence
	$$\lambda_{\phi_{2^s}}(i+2^s) = \abs{ \{k\ge 2^s \sep k\preceq i+2^s\} } = \abs{ \{k\ge 2^s \sep k\preceq i+2^{s+1}\} }$$
	where the last equality is implied by the fact that if $x$ is in $\{k\ge 2^s \sep k\preceq i+2^s\}$ then $x+2^s$ is in $\{k\ge 2^s \sep k\preceq i+2^{s+1}\}$ and vice versa. Then since $\lambda_{\phi_{2^s}}(i+2^{s+1})  = \abs{ \{k\ge 2^s \sep k\preceq i+2^{s+1}\} }$ and from the fact that $\lambda_{\phi_{2^s}}(i+2^{s+1}) = \lambda_{\phi_{2^s}}(i)$ for every $1<i<2^{s+1}$, we have
	$$\lambda_{\phi_{2^s}}(i)=\lambda_{\phi_{2^s}}(i+2^s)$$
	and hence $\lambda_{\phi_{2^s}}$ is eventually periodic of period $2^s$ starting from 1.\\
	Observe that $\lambda_{\phi_{2^s}}(i+2^s)=\abs{ \{k\ge 2^s \sep k\preceq i+2^s\} }$ is even for every $i$ different from $0$, therefore we have the formula $\phi_{2^s} = \sum_{j}\sigma_{j2^s}$.
\end{proof}

As a consequence of Theorem \ref{ThProd}, since $b(j^2s) \wedge b(i) = 0$ for each $i < 2^s$, we have $\sigma_{i+j2^s} = \sigma_i \cdot \sigma_{j2^s}$. We can use this to give a more concise formula for $\phi_t$, stated in the next corollary.
\begin{Cor}
	Let $2^{s-1}\le t\le 2^s-1$. Set $\psi_{2^s}=\phi_{2^s}+1$ and $\eta_t=\sum_{i=t}^{2^s}\sigma_i$. Then
	$$\phi_t=\psi_{2^s}\cdot \eta_t.$$
\end{Cor}

In the remaining part of this section we use $\phi_t$ to derive another family of symmetric polynomials. Let us consider the set of point of weight exactly $t$, and the polynomial $\rho_t$ vanishing at each point outside of this set.

\begin{Prop}\label{propRho}
	$\rho_t = \phi_{t+1}+ \phi_t = \phi_t \cdot (1 + \phi_{t+1})$.
\end{Prop}
\begin{proof}
	Apply the definition of $\phi_t$ and $\rho_t$ and use $\phi_{t_1} \cdot \phi_{t_2} = \phi_{\max(t_1 ,t_2 )}$.
\end{proof}

\begin{Th}
	$\rho_t$ is equal to the binary M{\"o}bius transform of $\sigma_{t}$.
\end{Th}
\begin{proof}
	The transform of $\sigma_t$ is exactly the polynomial vanishing at all points whose weight is different from $t$.
\end{proof}

We conclude this section with the following generalization of the idea
behind the derivation of $\rho_t$. A related result is shown in \cite{carlet2010boolean} to
characterise the Numerical Normal Form of binary symmetric functions.

\begin{Th}
	$\{\phi_t \}$ and $\{\rho_t \}$ are bases for the vector space of symmetric
	Boolean functions.
\end{Th}
\begin{proof}
	Since a symmetric Boolean function $f$ in $n$ variables assumes the same value on points whose Hamming weights are the same, it is determined by its simplified value vector $v_f$ of length $n+1$. We can write
	$$f = v_f(0) \rho_0 + v_f(1) \rho_1 + \dots + v_f(n)\rho_n.$$
	Since we can define $\rho_t$ in terms of $\phi_t$ and $\phi_{t-1}$ we also can write $f$ as a linear combination of $\phi_1,\dots, \phi_n$.
\end{proof}

\section{An application to linear codes}\label{Sect:Applic}

We recall that an $(n, M,d)_2$ nonlinear code is a subset of $\left(\mathbb{F}_2\right)^n$ of size $M$, where the Hamming distance between any pair of elements of $C$, called codewords, is at least equal to $d$. Any nonlinear code can be described as the image of a function $F : \left(\mathbb{F}_2 \right)^k \to \left(\mathbb{F}_2 \right)^n$ which we call the generator map of the code.  Without loss of generality we can assume that $F(0) = 0$, so that $0 \in C$. Indeed, if $F(0)=c_0\neq 0$ then we can consider instead of $C$ an equivalent code $\bar{C}$ obtained by the translation of $C$ by $c_0$, i.e. $\bar{c}\in\bar{C}$ if and only if $\bar{c}+c_0=c\in C$. The generator map of $\bar{C}$ is the map $\bar{F}=F+c_0$ and $\bar{F}(0)=0$. We remark that if $F$ is a linear map then $\mathrm{Im}(F)$ is a vector subspace of $\left(\mathbb{F}_2\right)^n$ and $C$ is thus a linear code. In this case we consider $k$ to be the dimension of $C$ as a vector subspace of $\left(\mathbb{F}_2\right)^n$ and thus $F$ is injective. Therefore, $F$ can be defined as a map of the form $m\mapsto m\cdot G$ where the $k\times n$ matrix $G$ is called generator matrix for $C$. The weight distribution of a code is the sequence of integers $A_0,\ldots,A_n$ with $A_i$ defined as the number of codewords of Hamming weight $i$. Observe that with our definition $A_0=1$, and for linear codes we have $A_i=0$ for any $1<i<d$, i.e. $A_i$ is the first non-zero element of the weight distribution after $A_0$. If we denote with $w$ the minimum non-zero Hamming weight of the codewords, then $w\leq d$ and in the linear case $w=d$.

We now show a way to determine the minimum weight of a code by using $\phi_t$. A related approach was proposed in \cite{guerrini2010computing} to the
systematic nonlinear case. The ideas behind the two methods are indeed
similar, even though the results of this section do not require the computation
of a Gröbner Basis.\\
We denote with $\phi_t^{(n)}$ the restriction of $\phi_t$ to the case of $n$ variables. 
\begin{Th}\label{ThCodes}
	Let $F:\left(\mathbb{F}_2\right)^k\to\left(\mathbb{F}_2\right)^n$ be the generator map of an $(n,M,d)_2$ code with minimum non-zero weight equal to $w$, i.e.
	$$
	\begin{array}{rccl}
		F:&\left(\mathbb{F}_2\right)^k&\to&\left(\mathbb{F}_2\right)^n\\
		&v&\mapsto & F(v)=(f_1(v),\ldots,f_n(v))\;,
	\end{array}
	$$
	with $F(0)=0$, $\left|\mathrm{Im}(F)\right|=M$ and
	$$
	\min_{c\in C\smallsetminus \{0\}}\{\mathrm{w}(c)\}=w\;.
	$$
	Then, $w\ge t $ if and only if $\phi_t^{(n)} \circ F = \phi_1^{(k)}$.
\end{Th}
\begin{proof}
	$w \ge t$ means that $\mathrm{w}(F(v))\ge t$ for each vector $v \neq 0 \in \left(\mathbb{F}_2\right)^k$, which can be written as
	$$\phi_t^{(n)}(F(v)) = 1 \text{ for each }v\neq 0 \in \left(\mathbb{F}_2\right)^k.$$
	This means that $\phi_t \circ F$ is the Boolean function in $k$ unknowns for which $0\mapsto 0$ and $0 \neq v \mapsto 1$, and this is exactly the definition of $\phi_1^{(k)}$.
\end{proof}
By considering the minimum distance instead of the minimum non-zero weight we derive the following corollary.
\begin{Cor}\label{CorCodes}
	Let $C$ be a nonlinear code as in Theorem \ref{ThCodes} and let $F$ be its generator map.\\
	Then, $d<t$ if $\phi_t^{(n)} \circ F \neq \phi_1^{(k)}$.
\end{Cor}
\begin{proof}
	If $\phi_t^{(n)} \circ F \neq \phi_1^{(k)}$, then by Theorem \ref{ThCodes} $w<t$, and by definition $d\le w$.
\end{proof}

We remark that $d$ can be strictly smaller than $w$. By applying Theorem \ref{ThCodes} at most $\log_2(n)$ times, we can determine precisely $w$, while Corollary \ref{CorCodes} states that we would only bound $d$. We can however restrict ourselves to linear codes, so that the minimum weight corresponds to the minimum distance. We obtain Corollary \ref{CorLinCodes}, whose proof is obtained by applying Theorem \ref{ThCodes} to linear codes.

\begin{Cor}\label{CorLinCodes}
	Let $G$ be the generator matrix of a linear $[n, k, d]_2$ code $C$ and let $F$ be the linear map
	$$
	\begin{array}{rccl}
		F:&\left(\mathbb{F}_2\right)^k&\to&\left(\mathbb{F}_2\right)^n\\
		&v&\mapsto & F(v)=v\cdot G\;.
	\end{array}
	$$
	Then, $d\ge t $ if and only if $\phi_t^{(n)} \circ F = \phi_1^{(k)}$.
\end{Cor}
We consider now the weight distribution of a linear code. Each codeword of weight $t$ is a root of $\rho_t+1$, therefore we have the following result.
\begin{Prop}
	Let $G$, $C$ and $F$ as in Corollary \ref{CorLinCodes}.\\
	The set of codewords of $C$ of weight $t$ correspond to the variety of the ideal generated by $1+\rho_t^{(n)}\circ F$.
\end{Prop}
From the above proposition we derive a result on the weight distribution of $C$.

\begin{Cor}
	Let $G$, $C$ and $F$ as in Corollary \ref{CorLinCodes}, and let $\{A_i\}$ be the weight distribution of $C$.\\
	Then, for any $i=1,\ldots,n$ the value $A_i$ is equal to the number of monomials in the ANF of the binary M{\"o}bius transform of $\rho_i^{(n)}\circ F$.
\end{Cor}
\begin{proof}
	Given a function $f$, the ANF of the binary M{\"o}bius transform of $f$ is the polynomial with coefficients given by the evaluation vector of $f$. Therefore, if $y=(y_1,\ldots,y_k)\in \left(\mathbb{F}_2\right)^k$ is such that $F(y)$ has weight $i$, then $y$ is a root of $1+\rho_i^{(n)}\circ F$. Then the function $\rho_i^{(n)}\circ F$ is equal to $1$ exactly on the vectors of $\left(\mathbb{F}_2\right)^k$ sent to codewords of $C$ of weight $i$. This implies that the monomial $X^y$ appears in the binary M{\"o}bius transform of $\rho_i^{(n)}\circ F$. By counting the monomials we conclude.
\end{proof}
\section{Conclusions}

From a theoretical point of view, the explicit description of $\phi_t$ allows the formulation of problems in which the solutions have requirements on their weight. Even though it is quite straightforward to simply check the weight of a given solution, in particular cases it could be an advantage to just add the right linear combination of $\rho_0,\dots, \rho_n$ to the generating system of the ideal.\\
Other than a theoretical overview of several properties of the polynomials $\phi_t$ and $\rho_t$, we have shown here how to obtain them either by applying the binary M{\"o}bius function or exploiting the periodicity of their sANF vectors. The latter allows to represent such polynomials using a minimum amount of memory. Finally, in Section \ref{Sect:Applic} we have shown an application of our results to Coding Theory, a novel theoretical method to check the minimum distance of a linear binary code.\\
We conclude by giving some remarks on our contribution to Coding Theory, even though the construction of dedicated algorithms and a study of the complexity of such procedures is beyond the purpose of this work. In the general case $[n, k]_2$ codes have no particular algebraic or geometric structure, i.e. the generator matrix of $C$ is chosen randomly. Then, since we are working with length $n$ codewords, the number of monomials in $\sigma_{i}$ is equal to $\binom{n}{i}$ and the computation of $\sigma_{i} \circ F$ (see Corollary \ref{CorLinCodes}) requires therefore $\binom{n}{i}$ multiplications, each one involving $i$ linear polynomials. This is however the worst case scenario. Dedicated algorithms could instead take advantage of the symmetric nature of $\phi_t$, and be designed to compute the minimum distance of particular classes of codes, for example for codes where the generator matrix present symmetries. A further development of this work can involve two directions:
\begin{itemize}
	\item the study of codes where the columns of their generator matrices present symmetries to be exploited;
	\item the study of polynomials, similar to $\phi_t$ and $\rho_t$, invariant under the action of non-trivial subgroups of the symmetric group. 
\end{itemize}

\subsection*{Acknowledgements}
	The first author acknowledges support from TIM S.p.A. through the PhD scholarship. The second author is a member of the INdAM Research group GNSAGA and of the Cryptography and Coding group of the Unione Matematica Italiana (UMI). The authors would like to thank Prof. Massimiliano Sala for the interesting hints on polynomial ideals and their applications to Coding Theory.

\bibliographystyle{plain}
\bibliography{vanishing_v1.bib}

\section*{Appendix}\label{Sect:Appendix}
We discuss here the structure of $\phi_t$ and $\rho_t$ for some particular $t$'s.

The Lucas' Theorem implies the following Lemma, useful to prove a particular property of $\phi_t$ for $t$ even.

\begin{Lem}\label{LemLuc}
	Let $i>j$ be two integers equal to $0 \mod 2^e$. Then $\binom{i}{j} \mod 2 =\binom{i+b}{j}$ for each $b<2^e$.
\end{Lem}
\begin{Prop}\label{Propmod2}
	Let $t=0 \mod 2^e$. Then if $\lambda_{\phi_t}$ is the sANF vector of $\phi_t$, then
	\begin{equation}\label{eqPro}
		\lambda_{\phi_t}(i)=0 \text{ for each } i\ne 0\mod 2^e.
	\end{equation}
\end{Prop}
\begin{proof}
	The case $i<t$ follows from Lemma \ref{Lems}, so let $i>t$ and let $\tilde{i} = j2^e$ be such that $\tilde{i}<i<\tilde{i}+2^e$.
	We will proceed by induction. Assume that Equation \eqref{eqPro}
	holds till the coefficient $\lambda_{\phi_t}(i-1)$. In this case we can write
	$$\phi_t = \sum_j \lambda_{\phi_t}(j2^e)\sigma_{j2^e}+\lambda_{\phi_t}(i)\sigma_{i} + R.$$
	We recall that since $R$ is symmetric, it is a linear combination of ESFs of degree larger than $i$.\\
	By definition, $\phi_t(v) = 1$ whenever $\mathrm{w}(v)\ge t$, so in particular we have both $\phi_t(v) = 1$ with $v$ of weight $i$ and $\phi_t (\tilde{v}) = 1$ with $\mathrm{w}(\tilde{v})=\tilde{i}$:
	\begin{equation*}
		\left\{
		\begin{array}{l}
			\phi_t(\tilde{v}) = \sum_j \lambda_{\phi_t}(j2^e)\sigma_{j2^e}(\tilde{v}) = 1 \\
			\phi_t(v) = \sum_j \lambda_{\phi_t}(j2^e)\sigma_{j2^e}(v) + \lambda_{\phi_t}(i)\sigma_i(v) = 1
		\end{array}
		\right.
	\end{equation*}
	By Lemma \ref{LemLuc} it follows that $\sigma_{j2^e}(v) = \sigma_{j2^e} (\tilde{v})$ for each $j$, so $\phi_t (v) = 1 + \lambda_{\phi_t}(i)\sigma_i(v) = 1$, hence $\lambda_{\phi_t}(i) = 0$.
\end{proof}

We summarise most of the results that we presented up to this point in the next theorems, which are methods to compute $\phi_t$ respectively for odd and even values of $t$.

\begin{Lem}\label{LemAppendix}
	Let $\mathrm{w}(v) = w$. Then $\sigma_i(v) = 1$ if and only if $\binom{w}{i} = 1 \mod 2$.
\end{Lem}
\begin{proof}
	Let $\{i_1 , \dots , i_w \}$ be the support of $v$, and let $m = x_{j_1} \cdots x_{j_i}$ be a monomial of degree $i$. Then $m(v) = 1$ if and only if $\{j_1 ,\dots , j_i \} \subseteq \{i_1 , \dots , i_w \}$. There are exactly $\binom{w}{i}$ such monomials of degree $i$, hence $\sigma_{i}(v) = \binom{w}{i} = 1 \mod 2$.
\end{proof}

\begin{Th}\label{ThOdd}
	Let $t$ be an odd integer, and let $s$ be such that $2^{s-1}<t\le 2^s$. Let $a_{t,t}=1$ and for $i\in\{ t+1,\dots, 2^s\}$ define
	\begin{equation}\label{EQai}
		a_{t,i} = 1 + \sum_{j=t}^{i-1}a_{t,j}\cdot \binom{i}{j} \mod 2.
	\end{equation}
	Then for $t\le i\le 2^s$ we have $\lambda_{\phi_t}(i)=a_{t,i}$, while for $0\le i < t$ we have $\lambda_{\phi_t}(i)=0$.
\end{Th}
\begin{proof}
	For $0\le i < t$, Lemma \ref{Lems} implies $\lambda_{\phi_t}(i)=0$. Equation \eqref{EQai} derives from the definition of $\phi_t$ and from Lemma \ref{LemAppendix} and the periodicity of $\lambda_{\phi_t}$ defines the whole sANF vector. 
\end{proof}

\begin{Th}
	Let $t = r2^e$ , with $r$ being an odd integer. Let $s$ be such that $2^{s-1} < r \le 2^s$. Let $\{a_{r,i} \}_i$ be the sequence of coefficients of $\phi_r$, as defined in \ref{ThOdd}. Then
	$$ \phi_t = \sum_{j=0}^{\infty}\sum_{i=t}^{2^s}a_{t,i}\sigma_{i+j2^e}. $$
\end{Th}

\begin{proof}
	We apply first Proposition \ref{Propmod2}, and then Theorem \ref{ThOdd}.
\end{proof}
From Theorem \ref{ThOdd} we can derive explicit formulas for some particular cases, related to $\phi_{2^s}$.
\begin{Cor}\label{Cort-}
	Let $t_{-} = 2^s - 1$ and let $t_+ = 2^{s-1} + 1$. Then
	\begin{align*}
		\phi_{t_-} &= \sum_j (\sigma_{2^s-1+j2^s} + \sigma_{2^s+j2^s})\\
		\phi_{t_+} &= \sum_j \sum_{i=t}^{2^s} \sigma_{i+j2^s}.
	\end{align*}
\end{Cor}

\begin{Cor}
	Let $t_- = 2^s-1$. Then
	\begin{align*}
		\rho_{2^s} &= \sum_j \sum_{i=2^s}^{2^{s+1}-1} \sigma_{i+j2^{s+1}}\\
		\rho_{t_-} &= \sum_j \sigma_{2^s-1+j2^s}.
	\end{align*}
\end{Cor}
\begin{proof}
	We apply Proposition \ref{propRho} to Corollary \ref{Cor2s} and Corollary \ref{Cort-}.
\end{proof}

\end{document}